\documentclass[a4paper,12 pt]{amsart}
\newtheorem{theo}{Theorem}[section]
\newtheorem{ex}[theo]{Example}

\newtheorem{fact}[theo]{Fact}
\newtheorem{question}[theo]{Question}
\newtheorem{prop}[theo]{Proposition}
\newtheorem{lem}[theo]{Lemma}
\newtheorem{cor}[theo]{Corollary}
\newtheorem{conj}[theo]{Conjecture}

\newtheorem{rema}[theo]{Remark}

\def \kbar {{\overline k}}
\def \Xbar {{\overline X}}

\def \Romannumeral #1 {\expandafter\uppercase\expandafter {\romannumeral #1} }

\def \pic {{\rm {Pic}}}
\def \div {{\rm{Div}}}
\def \Gal {{\rm{Gal}}}
\def \Alb {{\rm{Alb}}}

\def \calo {{\mathcal O}}
\def \spec {{\rm{Spec\,}}}

\def \dim {{\rm{dim\,}}}

\def \Hom {{\rm {Hom}}}
\def \Ext {{\rm {Ext}}}

\def \Aut{{\rm Aut \,}}

\def\ov{\overline}

\def \Z {{\bf Z}}
\def \Q {{\bf Q}}
\def \F {{\bf F}}

\def \NN {{\bf N}}

\def \Alb {{\rm Alb}}

\def \G {{\bf G}_m}

\def \Lim {\mathop{\lim}}
\def \dirlim {{\Lim\limits_\rightarrow\,}}
\def \invlim {{\Lim\limits_\leftarrow\,}}

\def\smallsquare{\vbox{\hrule\hbox{\vrule height 1 ex\kern 1 ex\vrule}\hrule}}
\def\enddem{\hfill \smallsquare\vskip 3mm}

\usepackage{amscd}
\usepackage{amssymb}
\usepackage{amsmath}
\usepackage{pb-diagram}

\everymath{\displaystyle}

\def \id {{\rm id}}

\DeclareFontFamily{U}{wncy}{}
\DeclareFontShape{U}{wncy}{m}{n}{%
   <5>wncyr5%
   <6>wncyr6%
   <7>wncyr7%
   <8>wncyr8%
   <9>wncyr9%
   <10>wncyr10%
   <11>wncyr10%
   <12>wncyr6%
   <14>wncyr7%
   <17>wncyr8%
   <20>wncyr10%
   <25>wncyr10}{}
\DeclareMathAlphabet{\cyrille}{U}{wncy}{m}{n}
\def\Sha{\cyrille X}

\title{Galois sections for abelianized fundamental groups}
\author{by David Harari and Tam\'as Szamuely}
\address{Universit\'e de Paris-Sud Math\'ematique, B\^atiment 425, 91405 Orsay, France}
\email{David.Harari@math.u-psud.fr}
\address{Alfr\'ed R\'enyi Institute of Mathematics, Hungarian Academy of Sciences, PO Box 127, H-1364 Budapest, Hungary}
\email{szamuely@renyi.hu}
\address{Mathematical Institute, University of Oxford, 24--29 St.\ Giles,
Oxford OX1 3LB, United Kingdom}
\email{flynn@maths.ox.ac.uk}\date{}
\begin{document}
\maketitle
\begin{center} {\em With an appendix by E. V. Flynn}
\end{center}\bigskip\bigskip

\section{Introduction}

Let $k$ be a perfect field with fixed algebraic closure $\kbar$, and $X$ a geometrically
integral separated scheme of finite type over $k$. Choose a geometric point $\bar
x:\spec\kbar\to X$. One then has the well-known exact sequence of profinite groups
\begin{equation}\label{piex}
1\to \pi_1(\Xbar, \bar x)\to \pi_1(X, \bar x)\to \Gal(\kbar|k)\to 1,
\end{equation}
where $\pi_1(X, \bar x)$ denotes the \'etale fundamental group  of $X$, and $\Xbar$ stands
for the base change $X\times_k\kbar$.

Each $k$-rational point $\spec\, k\to X$ determines a section of the structure map
$X\to\spec\, k$. As the \'etale fundamental group is functorial for morphisms of pointed
schemes, taking the induced map on fundamental groups defines a map
$$X(k)\to \{\text{\rm continuous sections of }  \pi_1(X, \bar x)\to \Gal(\kbar|k)\}/\sim,$$
where two sections are equivalent under the relation $\sim$ if they are conjugate under
the action of $\pi_1(\Xbar, \bar x)$.

According to the projective case of the famous Section Conjecture of Grothendieck
formulated in \cite{bf}, the above map should be a bijection when
 $k$ is a finitely generated field over $\Q$ and $X$ is a smooth projective curve of genus at least 2. Here injectivity is easy to prove and was known to Grothendieck; the hard and widely
open part is surjectivity.  Our main motivation for writing the present paper is the
following variant of Grothendieck's conjecture.

 \begin{conj}\label{conj} A  smooth projective curve $X$ of genus at least 2 over a number field $k$
 has a $k$-rational point if and only if the map $\pi_1(X, \bar x)\to \Gal(\kbar|k)$ has a continuous
section.
\end{conj}

As several people have reminded us, together with the theorem of Faltings ({\em quondam}
Mordell's Conjecture) this seemingly weaker statement actually implies the surjectivity
part of Grothendieck's Section Conjecture over a number field (see \cite{koenigsmann},
Lemma 1.7, itself based on an idea of Tamagawa \cite{tama}).

Oddly enough, until the recent preprint \cite{stix} of J. Stix there seems to have been no
 examples where the statement of Conjecture \ref{conj} holds in a nonobvious
way, i.e. where $X(k)=\emptyset$ and $\pi_1(X, \bar x)\to \Gal(\kbar|k)$ has no continuous
section. This may be because, as we shall see below, verifying that sections do {\em not}
exist is by no means straightforward. In the examples of Stix the curve $X$ has no points
over some non-archimedean completion of $k$. Thus the question remained whether there are
curves of genus at least 2 where $\pi_1(X, \bar x)\to \Gal(\kbar|k)$ has no continuous
section and at the same time $X$ has points everywhere locally, or in other words it is a
counterexample to the Hasse principle for rational points.

In this paper we answer this question via investigating sections of the abelianized
fundamental exact sequence of $X$. This is the short exact sequence
 \begin{equation}\label{piabex}
1\to \pi_1^{\rm ab}(\Xbar)\to \Pi^{\rm ab}(X)\to \Gal(\kbar|k)\to 1
\end{equation}
obtained from (\ref{piex}) by pushout via the natural surjection  $\pi_1(\Xbar, \bar
x)\to\pi_1^{\rm ab}(\Xbar)$, where $\pi_1^{\rm ab}(\Xbar)$ is the maximal abelian
profinite quotient of $\pi_1(\Xbar, \bar x)$. Of course, if (\ref{piex}) has a continuous
section, then so does (\ref{piabex}).

Many of our considerations will be valid in arbitrary dimension. In dimension greater than
one, however, a slight complication arises from torsion in the N\'eron-Severi group. To
deal with it, we introduce yet another pushout of exact sequence (\ref{piex}). Recall
(e.g. from \cite{kala}, Lemma 5) that for a smooth projective geometrically integral
$k$-scheme $X$  there is an exact sequence
\begin{equation}\label{pialb}
0\to S\to \pi_1^{\rm ab}(\Xbar)\to T({\rm Alb}_{X}(\kbar))\to 0
\end{equation}
of $\Gal(\kbar|k)$-modules, where ${\rm Alb}_{X}$ denotes the Albanese variety of $X$, the
notation $T({\rm Alb}_{X}(\kbar))$ means its full Tate module, and $S$ is the finite group
dual to the torsion subgroup of the N\'eron-Severi group $NS(\Xbar)$ of $\Xbar$.  Taking
the pushout of exact sequence (\ref{piabex}) by the surjection $\pi_1^{\rm
ab}(\Xbar)\twoheadrightarrow T({\rm Alb}_{X}(\kbar))$ we obtain an extension
  \begin{equation}\label{pi}
1\to T({\rm Alb}_{X}(\kbar))\to \Pi\to \Gal(\kbar|k)\to 1
\end{equation}
of profinite groups. In the case when $NS(\Xbar)$ is torsion-free (e.g. for curves or
abelian varieties) this sequence is the same as (\ref{piabex}).

Recall also (e.g. from \cite{es}, Theorem 4.1) that by the theory of Albanese varieties
over arbitrary perfect fields there is a canonical $k$-torsor ${\rm Alb}^1_X$ under ${\rm
Alb}_X$, characterized by the universal property that every morphism from $X$ to a
$k$-torsor under an abelian variety factors uniquely through ${\rm Alb}^1_X$. We may
finally state:

\begin{theo}\label{theodiv} Let $X$ be a smooth projective geometrically integral variety over a
perfect field $k$.

The map $\Pi\to \Gal(\kbar|k)$ has a continuous section if and only if the class of ${\rm
Alb}^1_X$ lies in the maximal divisible subgroup of the group $H^1(k,\Alb_X)$ of
isomorphism classes of $k$-torsors under ${\rm Alb}_X$.
\end{theo}

We remind the reader that the maximal divisible subgroup of an abelian group may be
smaller than the subgroup of infinitely divisible elements (the latter is not always a
divisible subgroup). Thus the condition of the theorem is stronger than just requiring
$\Alb_X^1$ to be infinitely divisible.

The condition of the theorem is related to zero-cycles on $X$. Namely, it follows from
Theorem 4.2 of \cite{es} that the {\em triviality} of the torsor $\Alb^1_X$ is equivalent
to the surjectivity of the map $CH_0(\Xbar)^{\Gal(\kbar|k)}\to \Z$, i.e. the existence of
a Galois-invariant zero-cycle class of degree 1 on $\Xbar$. In some cases the latter
condition is equivalent to the surjectivity of the map $CH_0(X)\to \Z$, i.e. the existence
of a zero-cycle of degree one on $X$. This is the situation, for instance, for curves over
a number field having points everywhere locally (see \cite{milne}, Proposition 2.5 or
\cite{es}, Proposition 3.2), or for curves of odd genus over a $p$-adic field
(\cite{licht}, Theorem 7 $c$).

We may thus view Theorem \ref{theodiv} as a link between the existence of a section for
(\ref{pi}) and the existence of a Galois-invariant zero-cycle class of degree 1 on
$\Xbar$. Over special fields we can say more about the relation of these two conditions.
When $k$ is a finite extension of $\Q_p$, we shall see that the two conditions are
equivalent for curves of genus $g$ with $p$ prime to $g-1$ (a fact essentially going back
to Lichtenbaum \cite{licht}). On the other hand, we shall give examples of other curves
and higher-dimensional varieties where (\ref{pi}) splits but there is no Galois-invariant
zero-cycle class. This yields the somewhat surprising fact that for curves over a $p$-adic
field the sections of exact sequence (\ref{piabex}) do not detect zero-cycles of degree 1
in general.

Over number fields we understand the situation less well but, inspired by work of
Bashmakov \cite{bash72}, we at least give examples where $H^1(k,\Alb_X)$ has trivial
maximal divisible subgroup, and therefore the splitting of (\ref{pi}) is equivalent to the
existence of a Galois-invariant zero-cycle class of degree 1 on $\Xbar$. In particular, we
shall point out the following instance of the failure of a local-global principle for
sections of exact sequence (\ref{piex}):

\begin{prop}\label{propflynn}
Let $X$ be a smooth projective curve over $\Q$ whose Jacobian is isogenous over $\Q$ to a
product of elliptic curves each of which has finite Tate-Shafarevich group and infinitely
many $\Q$-points. Assume moreover that $X$ has points everywhere locally but no
$\Q$-rational divisor class of degree 1. Then (\ref{piex}) has sections everywhere locally
but not globally.
\end{prop}

As Victor Flynn shows in the appendix, curves of genus $2$ satisfying the
assumptions of the proposition exist.\smallskip

In the last section we shall discuss the relation of the splitting of (\ref{pi}) to the
 elementary obstruction of Colliot-Th\'el\`ene and Sansuc for the existence of rational points on
varieties (\cite{bcts}, \cite{witt}).

We are indebted to Olivier Wittenberg for pointing out a mistake in a previous version of
this work and for several useful remarks. Jakob Stix, Michael Stoll and the referee also
made valuable comments. The second author acknowledges support from OTKA grant No. 61116
and from the BUDALGGEO Intra-European project. He also thanks the mathematics department
of the \'Ecole Normale Sup\'erieure, where part of this work was done, for its
hospitality.

\section{Proof of Theorem \ref{theodiv}}

In this section we prove Theorem \ref{theodiv}. Denote the Galois group $\Gal(\kbar|k)$ by
$\Gamma$. The starting observation is that the extension (\ref{pi}) of topological groups
gives rise to a class $[\Pi]$ in Tate's continuous cohomology group $H^2_{\rm
cont}(\Gamma, T(\Alb_X)(\kbar))$, defined by means of continuous cocycles. It is the
trivial class if and only if the projection $\Pi\to\Gamma$ has a continuous section.

There is a product decomposition
$$
H^2_{\rm cont}(\Gamma, T(\Alb_X)(\kbar))\cong \prod_{\ell} H^2_{\rm cont}(\Gamma,
T_\ell(\Alb_X)(\kbar))
$$
where $\ell$ runs over the set of all primes. Under this decomposition $[\Pi]$ corresponds
to a system of classes $[\Pi_\ell]\in H^2_{\rm cont}(\Gamma, T_\ell(\Alb_X)(\kbar))$, each
of which is the class of the extension
\begin{equation}\label{piell}
1\to T_\ell({\rm Alb}_{X}(\kbar))\to \Pi_\ell\to \Gamma\to 1
\end{equation}
obtained from (\ref{pi}) by pushout via the natural projection $T({\rm Alb}_{X}(\kbar))\to
T_\ell({\rm Alb}_{X}(\kbar))$. Now Theorem \ref{theodiv} follows from applying the
following proposition for each prime number $\ell$.

\begin{prop}\label{proptheodiv}
The class $[\Pi_\ell]\in H^2_{\rm cont}(\Gamma, T_\ell(\Alb_X)(\kbar))$ is trivial if and
only if the class of ${\rm Alb}^1_X$ lies in the maximal $\ell$-divisible subgroup of the
group $H^1(k,\Alb_X)$.
\end{prop}

For some of our considerations it will be more convenient to work with Jannsen's
continuous \'etale cohomology groups introduced in \cite{jannsen}. The coefficients in
Jannsen's theory are inverse systems $({\mathcal F}_n)$ of \'etale sheaves indexed by the
ordered set $\bf N$ of nonnegative integers. According to Theorem 2.2 of \cite{jannsen},
over $\spec\,k$ Jannsen's groups coincide with Tate's continuous Galois cohomology groups
for coefficient systems satisfying the Mittag-Leffler condition. This is the case for
systems of finite Galois modules, the only ones we shall consider.

It is useful to bear in mind the exact sequences
\begin{equation}\label{jannsenex}
0\to\invlim^1 H^{i-1}(k, {\mathcal F}_n) \to H^{i}_{\rm cont}(k, ({\mathcal F}_n))\to
\invlim H^{i} (k, {\mathcal F}_n)\to 0
\end{equation}
(see \cite{jannsen}, Proposition 1.6) for all $i>0$, where the left-hand-side group
involving the derived inverse limit is nontrivial in general.

The class $[\Pi_\ell]\in H^2_{\rm cont}(\Gamma, T_\ell(\Alb_X)(\kbar))$ we are interested
in can thus be viewed as a class in Jannsen's  group $H^2_{\rm cont}(k,
(_{\ell^n}\Alb_X))$; here ${}_{\ell^n}A$ denotes the $\ell^n$-torsion part of an abelian
group $A$. Denote by $(\Alb_X, \ell)$ the inverse system indexed by $\bf N$ where each
term is the \'etale sheaf $\Alb_X$ over $\spec\, k$ and the maps are given by
multiplication by $\ell$. There is an exact sequence of inverse systems of \'etale sheaves
\begin{equation}\label{exinv}
0\to (_{\ell^n}\Alb_X, \ell)\to (\Alb_X, \ell)\stackrel{\ell^n}\to (\Alb_X, \id)\to 0,
\end{equation}
where the last term is a constant inverse system. Here the map $\ell^n$ stands for the
isogeny given by multiplication by $\ell^n$; in particular, it is surjective. It yields a
coboundary map
\begin{equation}\label{deltaell}
\delta_\ell:\, H^1(k,\Alb_X)\to H^2_{\rm cont}(k, (_{\ell^n}\Alb_X))
\end{equation}
using the obvious fact that the continuous cohomology of the constant system
$(\Alb_X,\id)$ is just the usual cohomology of $\Alb_X$.

The key observation is now the following.

\begin{prop}\label{lemalb}
The map $\delta_\ell$ sends the class $[\Alb^1_X]\in H^1(k,\Alb_X)$ to $[\Pi_\ell]$.
\end{prop}

The proof of this fact uses a technical lemma from the literature that we copy here for
the readers' convenience.

\begin{lem} {\rm (\cite{skobook}, Lemma 2.4.5)} Let $1\to A\to B\to C\to 1$ be a central
extension of algebraic $k$-groups such that $B$, $C$ are geometrically connected and $A$
is finite. Let $Y$ be a right $k$-torsor under $C$. Choose a base point $\bar y_0\in
Y(\kbar)$, and let $\nu:\, \overline C\to \overline Y$ be the isomorphism of right torsors
under $\overline C$ sending the neutral element to $\bar y_0$. Let $\overline B\to
\overline Y$ be the composition of $\overline B\to\overline C$ with $\nu$. Then we have
the extension of groups
\begin{equation}\label{skorext} 1\to A(\kbar)=\Aut(\overline B|\overline Y)\to
\Aut(\overline B|Y)\to \Gamma\to 1 \end{equation}
 such that the induced $\Gamma=\Gal(\kbar|k)$-module structure on $A(\kbar)$ is its usual $\Gamma$-module
structure.

The class of the extension (\ref{skorext}) in $H^2(\Gamma,A)$ coincides with
$\partial([Y])$, where $\partial$ is the connecting homomorphism $H^1(\Gamma, C)\to
H^2(\Gamma,A)$.
\end{lem}

The proof given in \cite{skobook} actually shows more: the two classes coincide `at the
cocycle level'. More precisely, if we consider the 1-cocycle $c:\,\Gamma\to C(\kbar)$
representing $[Y]\in H^1(\Gamma, C)$ that comes from the identification of $\overline Y$
with $\overline C$ by means of $\bar y_0$, then in order to obtain a 2-cocycle
representing the image of $[Y]$ via $\partial$ one has to choose a continuous map $b:\,
\Gamma \to B(\kbar)$ lifting $c$. On the other hand, a 2-cocycle representing the class of
(\ref{skorext}) comes from the choice of a section of the projection $\Aut(\overline
B|Y)\to \Gamma$. The proof shows that the map $\gamma\mapsto (\bar x\mapsto b(\gamma)\cdot
\gamma(\bar x))$ is such a section, and the resulting 2-cocycles are the same.
\medskip

\noindent{\bf Proof of Proposition \ref{lemalb}:} Fix $n>0$, and denote by $\Pi_{\ell^n}$
the quotient of $\Pi_\ell$ obtained as the pushforward of (\ref{piell}) via the quotient
map $T_{\ell}(\Alb_X(\kbar))\to {}_{\ell^n}\Alb_X(\kbar)$. Denoting by $\lambda$ the
composite map $\Alb_\Xbar\stackrel{\ell^n}\to\Alb_\Xbar\to \Alb_X^1$, the extension
$$
1\to {}_{\ell^n}\Alb_X(\kbar)\to \Aut(\Alb_{\Xbar}\stackrel\lambda\to \Alb_X^1)\to
\Gamma\to 1
$$
identifies with
\begin{equation}\label{pil} 1\to {}_{\ell^n}\Alb_X(\kbar)\to
\Pi_{\ell^n}\to \Gamma\to 1.
\end{equation}
Indeed, this is so for $\Alb^1_X$ in place of $X$ by the description of the fundamental
group of an abelian variety, and then the Albanese map $X\to\Alb_X^1$ induces an
isomorphism of exact sequence (\ref{pil}) with the corresponding sequence for $\Alb_X^1$.

Now apply the lemma to the exact sequence
\begin{equation}\label{albext}
0\to {}_{\ell^n}\Alb_X\to\Alb_X\stackrel{\ell^n}\to\Alb_X\to 0 \end{equation}
 and the torsor $\Alb_X^1$ under $\Alb_X$. It says that the class of the extension (\ref{pil}) in $H^2(\Gamma,{}_{\ell^n}\Alb_X)$ is the
image of $[\Alb_X^1]$ by the coboundary map $H^1(\Gamma,\Alb_X)\to
H^2(\Gamma,{}_{\ell^n}\Alb_X(\kbar))$ coming from (\ref{albext}).

Finally, we exploit the remark made above that the identification just described holds at
the level of cocycles. Making $n$ vary we in fact obtain a coherent system of sections
$\Gamma\to\Pi_{\ell^n}$ coming from liftings of a 1-cocycle $\Gamma\to\Alb_X(\kbar)$
representing $[\Alb_X^1]$ to the $\ell^n$-coverings
$\Alb_\Xbar(\kbar)\stackrel{\ell^n}\to\Alb_\Xbar(\kbar)$. The 2-cocycles $\Gamma\times
\Gamma\to {}_{\ell^n}\Alb_X(\kbar)$ induced by these sections assemble to a continuous
2-cocycle $\Gamma\times \Gamma\to T_\ell(\Alb_X(\kbar))$ that describes the class of $\Pi$
in $H^2_{\rm cont}(\Gamma, T_\ell(\Alb_X(\kbar)))$. By construction, it also represents
the image of $[\Alb_X^1]$ by $\delta_\ell$.
\enddem

 Notice that in the last argument we had to work with cocycles instead of cohomology
groups because by exact sequence (\ref{jannsenex}) the group $H^2_{\rm cont}(\Gamma,
T_\ell(\Alb_X(\kbar)))$ itself is not necessarily the inverse limit of the
$H^2(\Gamma,{}_{\ell^n}\Alb_X(\kbar))$.
\medskip

\noindent{\bf Proof of Proposition \ref{proptheodiv}:} By the proposition just proven the
class $[\Pi_\ell]$ is trivial if and only if $[\Alb^1_X]$ lies in $\ker(\delta_\ell)$. It
thus suffices to show that $\ker(\delta_\ell)$ is the maximal $\ell$-divisible subgroup of
$H^1(k, \Alb_X)$. To see this, note first that since $H^2_{\rm cont}(k,
(_{\ell^n}\Alb_X))$ has no nonzero \mbox{$\ell$-divisible} subgroup (\cite{jannsen},
Corollary 4.9), the maximal $\ell$-divisible subgroup of $H^1(k, \Alb_X)$ lies in
$\ker(\delta_\ell)$.  On the other hand, the long exact sequence coming from (\ref{exinv})
shows that $\ker(\delta_\ell)$ equals the image of the map $H^1_{\rm cont}(k, (\Alb_X,
\ell))\to H^1(k, \Alb_X)$. We show that this image is an $\ell$-divisible group, which
will complete the proof.

There is a factorisation
$$
H^1_{\rm cont}(k, (\Alb_X, \ell))\to \invlim (H^1(k, \Alb_X), \ell) \to H^1(k, \Alb_X)
$$
where the middle term is the inverse limit of the inverse system given by multiplication
by $\ell$ on the group $H^1(k, \Alb_X)$. The first map comes from exact sequence
(\ref{jannsenex}), and it factors the map $H^1_{\rm cont}(k, (\Alb_X, \ell))\to H^1(k,
\Alb_X)$ because the inverse limit of the constant inverse system $(H^1(k, \Alb_X), \id)$
is of course $H^1(k, \Alb_X)$. Now the middle term is $\ell$-divisible by construction,
hence so is its image in $H^1(k, \Alb_X)$.
\enddem

\begin{rema}\rm Theorem \ref{theodiv} and its proof carry over {\em mutatis mutandis} to
the more general case when $X$ is only assumed to be smooth and quasi-projective. The
fundamental group has to be replaced by the tame fundamental group classifying covers
tamely ramified over the divisors at infinity, and  $\Alb_X^1$ has to be understood as
Serre's generalised Albanese torsor, which is universal for morphisms in torsors under
semi-abelian varieties over $k$ (the discussion in \cite{es}, Theorem 4.1 is in this
generality). The relation with the tame  fundamental group that generalises exact sequence
(\ref{pialb}) is explained, for instance, in \cite{ssz}, Proposition 4.4.\end{rema}

\section{The case of a $p$-adic base field}

Keeping the notations and assumptions of the previous section, recall that the {\em
period} of $X$ is defined as the order of the cokernel of the degree map
$CH_0(\Xbar)^\Gamma\to\Z$, and the {\em index} of $X$ as the order of the cokernel of the
degree map $CH_0(X)\to \Z$. If the index is 1, one says that $X$ has a {\em zero-cycle of
degree one}.

It is proven in Theorem 4.2 of \cite{es} that the period of $X$ equals the order of the
class $[\Alb^1_X]$ in $H^1(k, \Alb_X)$. Thus by Theorem \ref{theodiv} the surjectivity of
$CH_0(\Xbar)^\Gamma\to\Z$  implies that exact sequence (\ref{pi}) has a section. The
question arises whether the converse is true, or in other words, whether Galois-invariant
zero-cycles of degree 1 can be detected by sections of (\ref{pi}). By Theorem
\ref{theodiv} this amounts to asking:

\begin{question}\label{question}
Can $[\Alb^1_X]$  be a nonzero element in the maximal divisible subgroup of $H^1(k,
\Alb_X)$?
\end{question}

In this section we investigate the question in the case when $k$ is a finite extension of
$\Q_p$. The following well-known fact is crucial for our considerations:

\begin{fact}\label{mattuck}
If $A$ is an abelian variety over a $p$-adic field $k$, then $H^1(k, A)$ is isomorphic to
 $F\oplus (\Q_p/\Z_p)^r$ with $F$ finite and $r\geq 0$. \rm Indeed, by Tate's duality for
abelian varities over $p$-adic fields (\cite{adt}, Corollary I.3.4) the group $H^1(k,A)$
is the $\Q/\Z$-dual of $A^*(k)$, where $A^*$ denotes the dual abelian variety. By
Mattuck's theorem \cite{mattuck} the group $A^*(k)$ has a finite index open subgroup
isomorphic to $\Z_p^r$ for ${r=[k:\Q_p]\,\dim  A}$, whence the assertion.
\end{fact}

In particular, the group $H^1(k,\Alb_X)$ is the sum of a finite abelian group and a
$\Z_p$-module of finite cotype, the latter property meaning that the $\Q_p/\Z_p$-dual is a
finitely generated $\Z_p$-module. Therefore its maximal divisible subgroup coincides with
the subgroup of divisible elements, and the above question over $p$-adic $k$ is equivalent
to asking whether $[\Alb^1_X]$ can be a nonzero divisible element in $H^1(k, \Alb_X)$.

A first partial answer is contained in the next proposition. It is essentially a
reformulation of old results of Lichtenbaum \cite{licht}, and is already implicit in
\cite{stix}; the proof given here is different.

\begin{prop}
Let $k$ be a finite extension of $\Q_p$, and  $X$ a smooth proper geometrically connected
curve of genus $g \geq 1$. Assume that $p$ does not divide $(g-1)$. Then the map $\Pi^{\rm
ab}(X) \to \Gamma$ has a continuous section if and only if the degree map $CH_0(\ov
X)^{\Gamma} \to \Z$ is surjective. If $g$ is odd, this is equivalent to saying that $X$
has a zero-cycle of degree one.
\end{prop}

\begin{dem} By Fact \ref{mattuck} the group $H^1(k,\Alb_X)$ is isomorphic to $F \oplus (\Q_p/\Z_p)^r$ with $F$ finite.
Therefore a divisible class can only lie in the $(\Q_p/\Z_p)^r$-component, and as such
must be of order a power of $p$. By \cite{licht}, Theorem 7 $a$) the period of $X$, which
is the same as the order of the class $[\Alb_X ^1]$, divides $g-1$. Hence the assumption
that $p$ does not divide $(g-1)$ implies that the class $[\Alb_X ^1]$ is divisible if and
only if it is zero, i.e. when $X$ has period 1. Now the first assertion follows from
Theorem \ref{theodiv}. If moreover $g$ is odd, then according to \cite{licht}, Theorem 7
$c$) the fact that the period is $1$ entails that the index is also $1$, whence the second
assertion.
\end{dem}

We now turn to cases where the answer to Question \ref{question} is negative, i.e.  where
the class of $[\Alb^1_X]$ is nonzero and divisible in $H^1(k,\Alb_X)$. Such examples are
easy to construct: if $A$ is an abelian variety over the $p$-adic field $k$ such that
$A(k)$ is infinite, then the structure of $H^1(k, A)$ recalled in Fact \ref{mattuck}
implies that there is a nonzero divisible class in $H^1(k, A)$, coming from a torsor $X$
under $A$. By construction we have $\Alb_X=A$ and $\Alb^1_X=X$, and we are done.

The following proposition handles the more difficult task of constructing such elements
for curves of genus $>1$. Together with the previous discussion it shows that for curves
of genus $p+1$ the abelian quotient of the fundamental group does not always suffice to
detect zero-cycles of degree 1.

\begin{prop}
Let $p$ be an odd prime number. There exists a curve $Y$ of genus $p+1$ over $k=\Q_p$ such
that $\Alb^1_Y$ yields a nonzero divisible class in $H^1(k, \Alb_Y)$. Therefore the map
$\Pi^{\rm ab}(Y) \to \Gamma$ has a continuous section, but the degree map $CH_0(\ov
Y)^{\Gamma} \to \Z$ is not surjective.
\end{prop}

The proof below is inspired from a construction used by Sharif \cite{sharif}. His original
aim was to construct curves over $p$-adic fields with given index and period.\medskip

\begin{dem} To begin with, let $E$ be the Tate elliptic curve with parameter $p$ over $\Q_p$.
By the above argument we know that there is a nonzero divisible class in $H^1(\Q_p, E)$.
But we can be more specific. As $E(\Q_p)$ is isomorphic to $\Q_p ^*/p^{\Z} \cong \F_p^*
\oplus \Z_p$ by the multiplicative structure of $\Q_p$, we have $H^1(\Q_p,E)\cong
\Z/(p-1)\Z \oplus \Q_p/\Z_p$ by Tate's duality theorem recalled above. Consider a torsor
$X$ under $E$ whose class is a divisible element of order $p$ in $H^1(\Q_p,E)$.

We contend that $X$ is split by the degree $p$ cyclic and totally ramified extension $K$
of $\Q_p$.  To see this, write $K=\Q_p(\pi)$, with a uniformizing parameter $\pi$ of $K$
satisfying $N_{K/k}(\pi)=p$. The group $E(K)=K^*/p^{\Z}$ is isomorphic to  $\Z/p\Z \oplus
\F_p^* \oplus \Z_p$, and an explicit element can be written in the form $\pi^m
\zeta^i(1+pb)$, where $\zeta$ is a primitive $(p-1)$-st root of unity and $b\in \Z_p$. The
corestriction map $E(K) \to E(\Q_p)$ corresponds to multiplication by $p$ because the norm
of $\pi^m \zeta^i(1+pb)$ is $p^m \zeta^{pi}(1+pb)^p$ for every $m,i \in \Z$, $b \in \Z_p$.
Since the restriction map $H^1(\Q_p,E) \to H^1(K,E)$ is the dual of the corestriction map
$E(K) \to E(\Q_p)$ via Tate duality, we obtain that our order $p$ class $[X]$ is killed by
restriction to $H^1(K,E)$.

Let $f$ be a non-square element in the function field $\Q_p(X)$ of $X$, and let $Y$ be the
normalization of $X$ in $\Q_p(X)(\sqrt{f})$.  Consider the commutative diagram
$$
\begin{CD}
(\pic \, \ov X)^\Gamma @>{\rm deg}>> \Z @>>> H^1(\Q_p, \pic^0_X)\\
@VVV @VVV @VVV \\
(\pic \, \ov Y)^\Gamma @>{\rm deg}>> \Z @>>> H^1(\Q_p, \pic^0_Y)
\end{CD}
$$
coming from the exact sequence of $\Gamma$-modules
$$
0\to \pic^0_X(\overline\Q_p)\to \pic\, \ov X\to \Z\to 1
$$
and the similar one for $Y$. Since the pullback map $\pic \,\ov X \to \pic\, \ov Y$
multiplies degree by $2$, the middle vertical map in the above diagram is multiplication
by 2. Here $\pic^0_X=\Alb_X$ as we are dealing with curves, and similarly for $Y$. As
explained in the proof of \cite{es}, Theorem 4.2, the image of $1\in \Z$ in $H^1(\Q_p,
\pic^0_X)$ is the class of $\Alb_X^1$, and similarly for $Y$. We conclude that the map
$H^1(\Q_p, \Alb_X) \to H^1(\Q_p,\Alb_Y)$ sends $[X]=[\Alb_X ^1]$ to $2[\Alb_Y ^1]$. As
$[X]$ is a $p$-divisible element of order $p$ by construction, we obtain that if the order
of $[\Alb_Y ^1]$ is also $p$, then it must  be a nonzero divisible class in
$H^1(\Q_p,\Alb_Y)$ (recall that $p$ was assumed to be odd).

We finally show that it is possible to choose $f$ in such a way that the order of
$[\Alb^1_Y]$ is $p$. As in \cite{sharif}, for this it is sufficient to find an $f$ whose
number of ramification points is $2mp$ with $m$ odd (this yields a $Y$ of genus $mp+1$ by
the Hurwitz formula). To get an example with $m=1$, we proceed as in \cite{sharif}. Define
a divisor $D$ on $\ov X=\ov E$ by
$$D:=\sum_{\gamma \in \Gal(K|\Q_p)} \gamma(1)-\sum_{\gamma \in \Gal(K|\Q_p)}
\gamma(\pi)$$ where $\gamma(x)$ denotes the twisted action of $\Gamma$ on $\ov X(\overline
\Q_p)=E(\overline\Q_p)=\overline\Q_p^*/p^{\Z}$ by a cocycle $\xi$ inducing the class $[X]
\in H^1(\Gal(K|\Q_p),E)$ (recall that $[X]$ is split by $K$). The divisor $D$ is
Galois-equivariant on $\ov X$. As $N_{K/\Q_p}(\pi)=1/p$ is trivial in
$E(\Q_p)=\Q_p^*/p^{\Z}$, the same argument as in Lemma 8 of \cite{sharif} shows that $D$
is principal. Therefore we can find an $f$ such that ${\rm div}(f)=D$, which is the one we
were looking for.
\end{dem}

\section{Et resurrexit Bashmakov}

In this section we investigate Question \ref{question} in the case when $k$ is a number
field. It becomes especially interesting when $X$ has points over each completion of $k$,
for then the class $\Alb^1_X$ lies in the Tate--Shafarevich group $\Sha(\Alb_X)$, a group
that is conjecturally finite. In particular, a nontrivial class cannot be divisible in
$\Sha(\Alb_X)$. By Theorem \ref{theodiv}, if we knew that $\Sha(\Alb_X)$ intersects the
maximal divisible subgroup of $H^1(k, \Alb_X)$ trivially, we could conclude that exact
sequence (\ref{pi}) splits if and only if $\ov X$ has a Galois-invariant zero-cycle class
of degree 1. As remarked in the introduction, for $X$ of dimension 1 the latter condition
is equivalent to the existence of a degree 1 divisor on $X$.

But here we run into a serious problem: to our knowledge it is not known whether a nonzero
class in the Tate--Shafarevich group $\Sha(A)$ of an abelian variety $A$ over $k$ can be
divisible in the group $H^1(k, A)$. This phenomenon was studied by Bashmakov in his papers
\cite{bash64} and \cite{bash72} without deciding the issue either way. However, he gave
examples where $\Sha(A)$, if finite, intersects the maximal divisible subgroup of $H^1(k,
A)$ trivially. This is interesting from our point of view as it yields via Theorem
\ref{theodiv} examples where the map $\Pi\to\Gal(\kbar|k)$ has no section.

The following proposition is more or less implicit in the proof of (\cite{bash72}, Theorem
7). The approach here is different.

\begin{prop} \label{finiet} Let $A$ be an abelian variety over a number field $k$.
Fix an open subset $U=\spec(\calo_S)$ of the spectrum of the ring of integers in $k$ such
that $A$ and its dual $A^*$ extend to abelian schemes ${\mathcal A}$ and ${\mathcal A}^*$
over $U$, and let $\ell$ be a prime number invertible on $U$. Assume that the
Tate-Shafarevich group $\Sha(A)$ is finite. Then the following assertions are equivalent:

\smallskip

a) The $\ell$-primary torsion subgroup $H^1(U,{\mathcal A}^*)\{\ell\}$ is finite.

b) The closure of the image of the diagonal embedding $$A(k)\to\prod_{v \mid \ell}
A(k_v)$$ (for the product of $v$-adic topologies) is a subgroup of finite index.
\end{prop}

\begin{dem} Define $\Sha_\ell(A^*)$ as the subgroup of $H^1(k,A^*)$ consisting of those
elements whose restriction to $H^1(k_v,A^*)$ is zero for each $v \nmid\ell$. By
Cassels--Tate duality (\cite{adt}, Theorem I.6.13) the dual of the finite group $\Sha(A)$
is $\Sha(A^*)$, and by Tate's local duality (\cite{adt}, Corollary I.3.4) the dual of the
profinite group $H^0(k_v,A)$ is the discrete group $H^1(k_v,A^*)$ for each place $v$ (here
by convention for $v$ archimedean $H^0(k_v,A)$ means Tate's modified cohomology group).
Dualizing the exact sequence
$$0 \to \Sha(A^*) \to \Sha_\ell(A^*) \to \bigoplus_{v \mid\ell} H^1(k_v,A^*)$$
we therefore get an exact sequence
$$\prod_{v \mid\ell } H^0(k_v,A) \to
\Sha_\ell(A^*) ^D \to \Sha(A) \to 0$$ where the superscript $D$ means $\Q/\Z$-dual.

Now by the Cassels--Tate dual exact sequence (see e.g. \cite{hasza}, Prop. 5.3 for a more
general statement) this extends to an exact sequence
$$0 \to \overline{{A(k)}} \to \prod_{v \mid\ell} H^0(k_v,A) \to
\Sha_\ell(A^*)^D \to \Sha(A) \to 0$$ where $\overline{{A(k)}}$ denotes the closure of
$A(k)$ in the topological product. As $\Sha(A)$ was assumed to be finite, assertion $b)$
is equivalent to the finiteness of $\Sha_\ell(A^*)$.

\smallskip

On the other hand, by (\cite{adt}, Lemma II.5.5) there is an exact sequence
\begin{equation} \label{lram}
0 \to \Sha_\ell(A^*) \to H^1(U,{\mathcal A}^*) \to \bigoplus_{v \not \in U,\, v\nmid\ell}
H^1(k_v,A^*).
\end{equation}
It follows from Fact \ref{mattuck} that for $v$ a finite place not dividing $\ell$ the
$\ell$-primary torsion subgroup $H^1(k_v,A^*)\{\ell\}$ is finite, and this holds in the
archimedean case as well. Therefore by the sequence above the finiteness of
$\Sha_\ell(A^*)\{\ell\}$ is equivalent to the finiteness of $H^1(U,{\mathcal
A}^*)\{\ell\}$.  To finish the proof, it is sufficient to show that
$\Sha_\ell(A^*)\{\ell\}$ is of finite index in $\Sha_\ell(A^*)$. There is an exact
sequence
$$0 \to \Sha(A^*) \to \Sha_\ell(A^*) \to \bigoplus_{v \mid\ell}
H^1(k_v,A^*)$$ The group $\Sha(A^*)$ is finite (this follows from the finiteness of
$\Sha(A)$ by \cite{adt}, Remark I.6.14 $c)$), hence $\Sha(A^*)/\Sha(A^*)\{\ell\}$ is
finite and $\ell$-divisible. On the other hand, for $v$ dividing $\ell$ the group
$H^1(k_v,A^*)\{ \ell \}$ is of finite index in $H^1(k_v,A^*)$ (use again Fact
\ref{mattuck}), and the result follows.
\end{dem}

\begin{cor} If $E$ is an elliptic curve over $\Q$ with $\Sha(E)$ finite and $E(\Q)$
infinite, then $H^1(U, {\mathcal E})\{\ell\}$ is finite, with $U$ and $\ell$ as above and
$\mathcal E$ extending $E$ over $U$.
\end{cor}

\begin{dem} Since $E(\Q_\ell)$ has a finite index subgroup
isomorphic to $\Z_\ell$ by Fact \ref{mattuck}, the assumption that $E(\Q)$ is infinite
implies that the closure of $E(\Q)$ in $E(\Q_l)$ is of finite index. Now apply the
proposition.
\end{dem}

This is to be compared with Theorem 7 of \cite{bash72} where it is proven without assuming
the finiteness of $\Sha(E)$ that if $E(\Q)$ is infinite, an element of $\Sha(E)$ that
becomes infinitely $\ell$-divisible in $H^1(U, {\mathcal E})$ must be infinitely
$\ell$-divisible in $\Sha(E)$.\medskip

Based on the above corollary it is easy to construct higher dimensional examples.

\begin{cor} \label{isoelli}
Let $A$ be an abelian variety over $\Q$ that is isogenous over $\Q$ to a product of
elliptic curves $E_1 \times...\times E_r$ such that each $E_i$ has finite Tate-Shafarevich
group and infinitely many $\Q$-points. Then $H^1(U,{\mathcal A}^*)\{ \ell \}$ and
$H^1(U,{\mathcal A})\{ \ell \}$ are finite, with $U$ and $\ell$ as in the proposition.
\end{cor}

\begin{dem} The finiteness of both $\Sha(A)$ and of $H^1(U,{\mathcal A})\{ \ell \}$ is
preserved under isogeny, so we may assume $A=E_1 \times...\times E_r$. We then get product
decompositions of the groups $\Sha(A)$, $H^1(U,{\mathcal A})\{ \ell \}$ and
$H^1(U,{\mathcal A}^*)\{ \ell \}$ to which we may apply the previous corollary.
\end{dem}

\begin{prop} \label{pasdiv}
Let $A$ be an abelian variety over $\Q$ satisfying the assumptions of
Corollary~\ref{isoelli}. Then the maximal divisible subgroup of $H^1(\Q,A)$ is trivial.

Consequently, if $X$ is a smooth projective curve over $\Q$ whose Jacobian is of the above
type and which does not have a $\Q$-rational divisor class of degree 1, then the map
$\Pi^{\rm ab}(X)\to\Gal(\overline \Q|\Q)$ has no section.
\end{prop}

\begin{rema}\rm The proposition is especially interesting in the case when $X$ has a point (or at
least a degree 1 divisor) defined over each completion of $\Q$, for in this case the
projection $\Pi^{ab}(X)\to\Gal(\overline \Q|\Q)$ has a section everywhere locally but not
globally.

Michael Stoll has communicated to us the example of the hyperelliptic curve with affine
equation $y^2=3x^6 + 8x^4 + 2x^2 - 6$ whose Jacobian is the product of the rank 1 elliptic
curves $y^2 = x^3 - x^2 - 15x - 27$ and  $y^2 = x^3 - x^2 - 49x + 157$. This Jacobian is
odd in the terminology of \cite{poonen}, and therefore the curve has a zero-cycle of
degree one everywhere locally but not globally; it thus satisfies the above assumptions.

In the appendix Victor Flynn gives explicit examples of curves of genus $2$ that satisfy
the above assumptions and moreover have rational points everywhere locally. This then
entails that the projection $\pi_1(X,\bar x)\to\Gal(\overline \Q|\Q)$, and  not just the
abelianized variant, has a splitting everywhere locally but not globally. In Stoll's
example there is no point over $\Q_2$.

\end{rema}

For the proof of Proposition \ref{pasdiv} we need a general lemma.

\begin{lem} \label{udivisible}
Let $A$, $\mathcal A$ and $U$ be as in Proposition \ref{finiet}.  An element  $\alpha \in
H^1(U,{\mathcal A})$ of order invertible on $U$ is divisible in $H^1(U,{\mathcal A})$ if
and only if its image in $H^1(k,A)$ lies in the maximal divisible subgroup.
\end{lem}

\begin{dem} By decomposing $\alpha$ in its $\ell$-primary components we reduce to the case when
$\alpha$ has $\ell$-power order for a prime $\ell$ invertible on $U$. By \cite{adt}, Lemma
II.5.5 the $\ell$-primary component $H^1(U,{\mathcal A}) \{\ell\}$ of the torsion group
$H^1(U,{\mathcal A})$ is of finite cotype, so its maximal divisible subgroup coincides
with the subgroup of divisible elements. As the image of a divisible subgroup by a
homomorphism of abelian groups is again divisible, the necessity of the condition of the
proposition follows. For sufficiency recall again from (\cite{adt}, Lemma II.5.5) the
exact sequence
$$
0 \to H^1(U,{\mathcal A}) \to H^1(k,A) \to \bigoplus_{v \in U} H^1(k_v,A).
$$
Here for $v \in U$ the group $H^1(k_v,A)\{\ell\}$ is finite by Fact \ref{mattuck} and our
assumption that $v$ does not divide $\ell$. Therefore the last term of the sequence
contains no nonzero $\ell$-divisible subgroup, so the maximal divisible subgroup of
$H^1(k,A)\{\ell\}$ is contained in $H^1(U,{\mathcal A})$. This is what we wanted to show.
\end{dem}

\noindent{\em Proof of Proposition \ref{pasdiv}.} The second assertion follows from the
first via Theorem \ref{theodiv} and the fact, already noted in the introduction, that $X$
has a $\Q$-rational divisor class of degree 1 if and only if $\Alb^1_X$ has trivial class
in $H^1(\Q,A)$.

For the first assertion pick $\alpha$ from the maximal divisible subgroup of $H^1(\Q,A)$.
We may find an open subset $U \subset \spec(\calo_k)$ such that $A$ extends to an abelian
scheme ${\mathcal A}$ over $U$ and $\alpha$ extends to an element $\alpha_U\in
H^1(U,\mathcal A)$ of order invertible on $U$. By the lemma $\alpha_U$ is divisible in
$H^1(U,\mathcal A)$. But by Corollary~\ref{isoelli} the group $H^1(U,{\mathcal A})\{ \ell
\}$ is finite for each prime $\ell$ dividing the order of $\alpha_U$, so $\alpha_U=0$ as
required.
\enddem

\section{Relation with the elementary obstruction}

In this section we investigate a different kind of condition for the vanishing of the
class $[\Pi]$. We begin with a continuous analogue for the well-known Ext spectral
sequence in usual Galois cohomology.

\begin{prop}
Let $k$ be a field, $\Gamma=\Gal(k_s|k)$,  $(M_n)$ a direct system of $\Gamma$-modules
indexed by $\NN$, and $N$ a $\Gamma$-module. There is a spectral sequence
\begin{equation}\label{ss}
H^p_{\rm cont}(k, (Ext^q(M_n, N))\Rightarrow \Ext^{p+q}_k(\dirlim M_n, N).
\end{equation}
\end{prop}

 Here the notation  $Ext^q(A, B)$ stands for the $q$-th left derived
functor of the functor $A\mapsto Hom(A,B)$, where $Hom(A, B)$ consists of those
homomorphisms $A\to B$ that are invariant under some open subgroup of $\Gamma$. Under the
correspondence between $\Gamma$-modules and \'etale sheaves on $\spec\, k$ this functor
corresponds to the inner Hom in the category of \'etale sheaves on $\spec\, k$. \smallskip

\begin{dem} This is the continuous analogue of the spectral sequence of \cite{adt}, Example 0.8.
As in the discrete case, it arises as a spectral sequence of composite functors. Let $F$
be the functor from the category of  $\Gamma$-modules to the category of inverse systems
of $\Gamma$-modules given by
$$
F:\,N \mapsto (Hom(M_n, N)),
$$
and $G$ the functor from the category of inverse systems of $\Gamma$-modules to that of
abelian groups given by
$$
G:\, (P_n)\mapsto (\invlim P_n)^\Gamma.
$$
According to Jannsen's definition in \cite{jannsen}, the $q$-th right derived functor of
$G$ maps $(P_n)$ to $H^q_{\rm cont}(k, (P_n))$.

We now check that $F$ transforms injective objects to $G$-acyclic objects. If $N$ is
injective, then for each $n$ the $\Gamma$-module $Hom(M_n, N)$ is acyclic for the functor
$A\to A^\Gamma$, as we see by applying Lemma 0.6 of \cite{adt} with $M=\Z$, $H=\{1\}$ and
$M_n$ and $N$ in place of $N$ and $I$. Therefore by \cite{jannsen}, Proposition 1.2 the
system $(Hom(M_n,N))$ is acyclic for $G$.

By the previous paragraph a spectral sequence for the composite functor $G\circ F$ exists.
Its $E_2^{pq}$-term is $H^p_{\rm cont}(k, (Ext^q(M_n, N)))$, and its abutment is the
$(p+q)$-th right derived functor of $$G\circ F:\,N\mapsto (\invlim Hom(M_n, N))^\Gamma.$$
But
\begin{align*}
(\invlim Hom(M_n, N))^\Gamma &\cong\invlim (Hom(M_n, N))^\Gamma\cong\invlim
(\Hom_\kbar(M_n, N))^\Gamma\\ \cong(\invlim \Hom_\kbar(M_n, N))^\Gamma
&\cong(\Hom_\kbar(\dirlim M_n, N))^\Gamma\cong\Hom_k(\dirlim M_n, N).
\end{align*}

\end{dem}

\begin{cor}\label{corext}
Let $k$ be a field, $\Gamma=\Gal(k_s|k)$ and $(M_n)$ a direct system of {\em finitely
generated} $\Gamma$-modules. For all divisible $\Gamma$-modules $N$ there are isomorphisms
$$
H^i_{\rm cont}(k, (\Hom_\kbar(M_n, N)))\cong \Ext^i_k(\dirlim M_n, N).
$$
\end{cor}

\begin{dem} For arbitrary $\Gamma$-modules $N$ we have $Ext^q(M_n, N)=\Ext^q_\kbar(M_n, N)$ for all $q$ and $n$,
because this is so for $q=0$ when $M_n$ is finitely generated.  If moreover $N$ is
divisible, the latter groups are trivial for $q>0$.

\end{dem}

We now apply the above corollary in a concrete situation. {\em Until the end of the
section we assume the base field $k$ is of characteristic 0.}

\begin{cor}\label{homlem} Let $X$ be a smooth projective geometrically integral variety over $k$.
There exists a canonical isomorphism of abelian groups
$$
H^2_{\rm cont}(k, T({\rm Alb}_{X}(\kbar))\cong \Ext^2_k((\pic^0\,\Xbar)_{\rm tors},
\kbar^\times)
$$
the extension group being taken in the category of $\Gamma=\Gal(\kbar|k)$-modules.
\end{cor}

Here $A_{\rm tors}$ denotes the torsion subgroup of an abelian group $A$.\medskip

\begin{dem} Fix a prime number $\ell$, and apply Corollary \ref{corext} with $i=2$,
${M_n={}_{\ell^n}{\pic}^0_X(\kbar)}$ (the $\ell^n$-torsion subgroup of
${\pic}^0_X(\kbar)$) and $N=\kbar^\times$. Then use the Weil pairing between the
$\ell^n$-torsion of Albanese and Picard varieties to obtain
$$
H^2_{\rm cont}(k, T_\ell({\rm Alb}_{X}(\kbar)))\cong \Ext^2_k((\pic^0\,\Xbar)_{\ell-\rm
tors}, \kbar^\times).
$$
The corollary follows by taking direct products over all $\ell$.
\end{dem}

Now recall (e.g. from \cite{skobook}, (2.16)) that there is a spectral sequence
$$
E_2^{pq}=\Ext^p_k((\pic^0\,\Xbar)_{\rm tors}, H^q(\Xbar, \G))\Rightarrow
\Ext^{p+q}_X(\pi^*({\pic^0_X})_{\rm tors}, \G),
$$
where $\pi:\, X\to\spec\, k$ is the natural projection and the abutment is an extension
group of \'etale sheaves on $X$. The differential $d:\, E_2^{01}\to E_2^{20}$ translates
as a map
$$
d:\, \Hom_\Gamma((\pic^0\,\Xbar)_{\rm tors}, \pic\,\Xbar)\to \Ext^2_k((\pic^0\,\Xbar)_{\rm
tors}, \kbar^\times).
$$
The natural inclusion map $i:\, (\pic^0\,\Xbar)_{\rm tors}\to \pic\,\Xbar$ thus yields a
class $d(i)\in\Ext^2_k((\pic^0\,\Xbar)_{\rm tors}, \kbar^\times)$.

\begin{prop}\label{piclass} Under the isomorphism of Corollary \ref{homlem} the class $[\Pi]$ in $H^2_{\rm cont}(k,
T({\rm Alb}_{X}(\kbar))$ corresponds to the class $d(i)$ in the group
$\Ext^2_k((\pic^0\,\Xbar)_{\rm tors}, \kbar^\times)$.
\end{prop}

\begin{dem} According to (\cite{skobook}, Theorem 2.3.4 $a)$), the map
$$
d':\,\Hom_\Gamma(\pic\,\Xbar, \pic\,\Xbar)\to \Ext^2_k(\pic\,\Xbar, \kbar^\times)
$$
coming from the spectral sequence
$$
\Ext^p_k(\pic\,\Xbar, H^q(\Xbar, \G))\Rightarrow \Ext^{p+q}_X(\pi^*{\pic\,X}, \G)
$$
sends the identity map of $\pic\, \Xbar$ to the opposite of the class of the 2-extension
\begin{equation}\label{ex}
1\to\kbar^\times\to \kbar(\Xbar)^\times\to \div\,\Xbar\to \pic \,\Xbar\to 1.
\end{equation}
The inclusion of $\Gamma$-modules $(\pic^0\,\Xbar)_{\rm tors} \to \pic\, \Xbar$ gives rise
to a commutative diagram
$$
\begin{CD}
\Hom_\Gamma(\pic\,\Xbar, \pic\,\Xbar) @>{d'}>> \Ext^2_k(\pic\,\Xbar, \kbar^\times)\\
@VVV @VVV \\
\Hom_\Gamma((\pic^0\,\Xbar)_{\rm tors}, \pic\,\Xbar) @>d>> \Ext^2_k((\pic^0\,\Xbar)_{\rm
tors}, \kbar^\times).
\end{CD}
$$
Therefore $-d(i)$ is the image of the class of (\ref{ex}) by the right vertical map. The
map in question factors as
$$
\Ext^2_k(\pic\,\Xbar, \kbar^\times)\to \Ext^2_k(\pic^0\,\Xbar,
\kbar^\times)\stackrel\lambda\to \Ext^2_k((\pic^0\,\Xbar)_{\rm tors}, \kbar^\times),
$$
and there is a natural map $\Phi:\,\Ext^2_{k_{\text{\rm \'et}}}(\pic^0_X, \G)\to
\Ext^2_k(\pic^0\,\Xbar, \kbar^\times)$, where the first group is an Ext-group in the
category of sheaves on the big \'etale site of $k$, and $\Phi$ is induced by restriction
of sheaves to the small \'etale site. The 2-extension (\ref{ex}) extends to a natural
2-extension
$$1\to\G\to {\mathcal K}_X^\times\to \div_X\to \pic_X\to 1$$
 of sheaves on the big \'etale site of $k$. Its pullback
via the map $\pic^0_X\to \pic_X$ yields a class $c\in\Ext^2_{k_{\text{\rm
\'et}}}(\pic^0_X, \G)$ whose image by $\lambda\circ\Phi$ is $-d(i)$.

Now recall that there is a natural isomorphism
$$
\phi:\,\Ext^2_{k_{\text{\rm\'et}}}(\pic^0_X, \G)\stackrel\sim\to H^1(k, {\rm Alb}_{X})
$$
coming from the local-to-global spectral sequence
\begin{equation}\label{lg}
H^p(k, Ext^i_{k_{\text{\rm \'et}}}(\pic^0_X, \G))\Rightarrow
\Ext^{p+q}_{k_{\text{\rm\'et}}}(\pic^0_X, \G),
\end{equation}
the Barsotti-Weil isomorphism $Ext^1_{k_{\text{\rm \'et}}}(\pic^0_X, \G)\cong {\rm
Alb}_{X}$ and the vanishing of $Ext^i_{k_{\text{\rm \'et}}}(\pic^0_X, \G)$ for $i\neq 1$
(see \cite{hasza}, Remark 4.1 for discussion and references on this point).

Skorobogatov checked in (\cite{skoro}, Proposition 2.1) that $\phi(c)$ is the opposite of
the class of $\Alb^1_X$. Now consider the diagram
\begin{equation}\label{keydiag}
\begin{CD}
\Ext^2_{k_{\text{\rm\'et}}}(\pic^0_X, \G) @>\phi>> H^1(k, {\rm Alb}_{X})\\
@V{\Lambda\circ\Phi}VV @VV{\delta}V \\
\Ext^2_k((\pic^0\,\Xbar)_{\rm tors}, \kbar^\times) @>>> H^2_{\rm cont}(k,
T(\Alb_X(\kbar)))
\end{CD}
\end{equation}
where $\delta$ is the product of the maps $\delta_\ell$ introduced in (\ref{deltaell}).
Once we have proven that the diagram commutes, it will follow from the above
considerations and Proposition \ref{lemalb} that $d(i)=-(\lambda\circ\Phi)(c)$ maps to
$[\Pi]$ in $H^2_{\rm cont}(k, T(\Alb_X(\kbar))).$

This commutativity is, however, a formal exercise in derived categories. As before, we
shall check it by writing $\Ext^2_k((\pic^0\,\Xbar)_{\rm tors}, \kbar^\times)$ as the
product of its $\ell$-primary components $\Ext^2_k((\pic^0\,\Xbar)_{\ell-{\rm tors}},
\kbar^\times)$, $T(\Alb_X(\kbar)$ as the product of the $T_\ell(\Alb_X(\kbar)$, and
considering each $\delta_\ell$ separately. To ease notation, we set $M=\pic^0_X$, $N=\G$,
$M_n={}_{\ell^n}\pic^0_X$, these being sheaves on the big \'etale site
$k_{\text{\rm\'et}}$.

The upper horizontal map $\phi$ then comes from the isomorphism
$$
{\bf R}\Hom_{k_{\text{\rm\'et}}}(M, N)\cong {\bf R}\Gamma(k,{\bf
R}Hom_{k_{\text{\rm\'et}}}(M, N)),
$$
which is the derived category version of spectral sequence (\ref{lg}). Similarly, the
lower horizontal map comes from
$$
{\bf R}\Hom_{k}(\dirlim M_n, N)\cong {\bf R}(\invlim{\bf R}Hom_{k}(M_n, N)^\Gamma),
$$
the derived category version of spectral sequence (\ref{ss}). Here the $\Hom$- and
$Hom$-functors concern the restrictions of sheaves to the small \'etale site of $k$. The
natural map $\dirlim M_n \to M$ induces a commutative diagram
$$
{\bf R}\Hom_{k_{\text{\rm\'et}}}(M, N)\cong {\bf R}\Gamma(k,{\bf
R}Hom_{k_{\text{\rm\'et}}}(M, N))
$$
$$
\downarrow\qquad\qquad\qquad\qquad\quad\downarrow
$$
$$
{\bf R}\Hom_{k}(\dirlim M_n, N)\cong {\bf R}(\invlim{\bf R}Hom_{k}(M_n, N)^\Gamma).
$$
We claim that the $\ell$-primary part of diagram (\ref{keydiag}) arises from the above by
applying the functor $H^2$. Indeed, we have already recalled above the vanishing of
$Ext^i_{k_{\text{\rm\'et}}}(M,N)$ for $i\neq 1$, so
$$H^j(k,{\bf
R}Hom_{k_{\text{\rm\'et}}}(M, N))\cong H^{j-1}(k, Ext^1_{k_{\text{\rm\'et}}}(M, N))$$ for
all $j>0$. The Barsotti--Weil formula then identifies the right-hand-side group for $j=2$
with the upper right group in (\ref{keydiag}). Moreover, we have seen in the proof of
Corollary \ref{corext} that the groups $Ext^i(M_n, N)$ vanish for $i>0$, whence
isomorphisms
$$
H^j_{\rm cont}(k, ({\bf R}Hom_{k}(M_n, N)))\cong H^j_{\rm cont}(k, (Hom_{k}(M_n, N)))
$$
and we obtain the lower right group in (\ref{keydiag}) for $j=2$. The only point that
still needs some justification is the identification of the map
$$
H^{1}(k, Ext^1_{k_{\text{\rm\'et}}}(\pic^0_X, \G))\to H^2_{\rm cont}
(k,(\Hom_k({}_{\ell^n}\pic^0_X, \kbar^\times)))
$$
with $\delta_\ell$. This is because the Barsotti--Weil formula translates the exact
sequence
$$
0\to {}_{\ell^n}\Alb_X\to \Alb_X\stackrel{\ell^n}\to\Alb_X\to 0
$$
of sheaves on $k_{\text{\rm\'et}}$ to
$$
0\to Hom_{k_{\text{\rm\'et}}}({}_{\ell^n}\pic^0_X, \G)\to
Ext^1_{k_{\text{\rm\'et}}}(\pic^0_X, \G)\stackrel{\ell^n}\to
Ext^1_{k_{\text{\rm\'et}}}(\pic^0_X, \G)\to 0
$$
in view of $Hom_{k_{\text{\rm\'et}}}(\pic^0_X, \G)=0$.
\end{dem}

According to  (\cite{skobook}, Theorem 2.3.4 $b)$) the class (\ref{ex}) is trivial if and
only if the exact sequence
\begin{equation}\label{obx}
1\to\kbar^\times\to \kbar(\Xbar)^\times\to \kbar(\Xbar)^\times/\kbar^\times\to 1
\end{equation}
has a Galois-equivariant section. The latter property is called, following
Colliot-Th\'el\`ene and Sansuc, the vanishing of the {\em elementary obstruction}.
Therefore Proposition \ref{piclass} implies:

\begin{cor} If the elementary obstruction vanishes,
the projection $\Pi\to\Gal(\kbar|k)$ has a section.
\end{cor}

Note that the converse is not true (there are plenty of rational, hence geometrically
simply connected varieties with nontrivial elementary obstruction).

\begin{rema}\label{remeo}\rm ${}$\smallskip
An argument similar to the one proving Proposition \ref{homlem} gives an isomorphism
$$
H^2_{\rm cont}(k, \Pi^{ab}(X))\cong \Ext^2_k((\pic\, \Xbar)_{\rm tors}, \kbar^\times).
$$
Whence reasoning as above we obtain that the vanishing of the elementary obstruction
implies the existence of a section for the projection $\Pi^{\rm ab}(X)\to\Gal(\kbar|k)$ as
well. H\'el\`ene Esnault and Olivier Wittenberg tell us that they have independently
obtained a similar result.
\end{rema}

\section{Appendix: An example satisfying the conditions of
Proposition 1.3}

\begin{center} {\em by E. V. Flynn}
\end{center}\bigskip

In this appendix, we shall show the existence of an example
which satisfies the conditions of Proposition \ref{propflynn}, that is,
a smooth projective curve over~$\mathbb Q$, with points
everywhere locally, with no $\mathbb Q$-rational divisor
class of degree~1, and whose Jacobian is isogenous over~$\mathbb Q$
to a product of elliptic curves (each defined over~$\mathbb Q$),
each of which has finite Tate-Shafarevich group and infinitely
many $\mathbb Q$-rational points.
We first recall the following theorem of Kolyvagin~\cite{koly:anal}.
\begin{theo}\label{thm:koly}
Let~$\mathcal{E}$ be an elliptic curve, defined over~$\mathbb Q$.
If~${\mathcal{E}}({\mathbb Q})$ has analytic rank~0 or~1, then
$\Sha({\mathcal{E}}/{\mathbb Q})$ is finite and the rank of ${\mathcal{E}}({\mathbb Q})$
is equal to its analytic rank.
\end{theo}
\par
The above result was originally proved in~\cite{koly:anal},
subject to the condition that~${\mathcal{E}}$ is modular;
since we now know~\cite{bcdt:mod} that all elliptic curves
over~$\mathbb Q$ are modular, the result can be given,
as above, without needing to state the modularity condition.
It is then natural to seek an example for which the Jacobian
is isogenous over~$\mathbb Q$ to a product of elliptic curves,
each of rank~1. There are various examples in the
literature (for example,
%\cite{blm:manin},
%\cite{cass:q},
\cite{flynn:manin})
%\cite{ps:ct},
%\cite{sk:shim}
of curves violating the Hasse principle, which make use of the
existence or non-existence of a rational divisor class of
degree~1; however, none of these examples satisfy the conditions
of Proposition~1.3. It seems likely that some of the curves
in~\cite{bruinstoll:manin} should satisfy these conditions,
however it is difficult to identify from the published data
which specific curves do so.
\par
We shall seek a genus~2 example
\begin{equation}\label{eq:genus2}
{\mathcal C} : Y^2 = F(X) = f_6 X^6 + f_5 X^5 + f_4 X^4 + f_3 X^3
+ f_2 X^2 + f_1 X + f_0,
\end{equation}
where $f_0,\ldots ,f_6\in {\mathbb Q}$ and $f_6\not= 0$.
Let~$J$ denote the Jacobian of~$\mathcal{C}$.
Following Chapter~1 of~\cite{prolegom},
the model~(\ref{eq:genus2}) should be taken as shorthand
for the corresponding projective smooth curve,
and we let~$\infty^+, \infty^-$ denote
the points on the non-singular curve that lie over the
singular point on~(\ref{eq:genus2}) at infinity.
For any field~$k$, with ${\mathbb Q} \subseteq k$,
we have~$\infty^+, \infty^- \in {\mathcal{C}}(k)$
when the coefficient of~$X^6$ is a square in~$k$.
We shall adopt the customary
shorthand notation $\{ P_1, P_2 \}$ to denote
the divisor class $[ P_1 + P_2 - \infty^+ - \infty^-]$,
which is in $J({\mathbb Q})$ when $P_1,P_2$ are points
on~$\mathcal{C}$ and either $P_1,P_2$ are both $\mathbb Q$-rational or
$P_1,P_2$ are quadratic over~$\mathbb Q$ and conjugate.
Let $F(X) = F_1(X)\ldots F_m(X)$ be the factorisation
of~$F(X)$ into irreducible polynomials over~$\mathbb Q$;
for each~$i$, let $\theta_i$ be a root of~$F_i(X)$
and let $L_i = {\mathbb Q}(\theta_i)$. Following p.49 of~\cite{prolegom},
we define the homomorphism
\begin{equation} \label{eq:mumapgen2}
\begin{array}{rl}
\mu &: J({\mathbb Q}) \rightarrow
   \Bigl( L_1^*/ (L_1^*)^2\times
           \ldots \times L_m^*/ (L_m^*)^2 \Bigr) /
{\lower2.5pt\hbox{$\widetilde{\ \ }$}},\\
& \{(x_1,y_1),(x_2,y_2)\} \mapsto [ (x_1 - \theta_1)(x_2 - \theta_1), \ldots ,
(x_1 - \theta_m)(x_2 - \theta_m) ],\\
\end{array}
\end{equation}
\noindent
where the equivalence relation ${\lower2.5pt\hbox{$\widetilde{\ \ }$}}$
is defined by
$$
[ a_1 ,\ldots , a_m ] {\lower2.5pt\hbox{$\widetilde{\ \ }$}}
[b_1 ,\ldots , b_m ] \Leftrightarrow
a_1 = w b_1 , \ldots , a_m = w b_m, \hbox{ for some }
w\in {\mathbb Q}^*.
$$
Since~$\mu$ is a map to a Boolean group, its
kernel clearly contains~$2J({\mathbb Q})$.
The following result (Lemma~6 in~\cite{flynn:manin}) gives a
way of showing the non-existence of a
$\mathbb Q$-rational divisor class of degree~1.
\begin{lem}\label{lem:SMALLKERMUimplyNORDCDO}
Let $\mathcal{C} : Y^2 = F(X)$ be as
in~(\ref{eq:genus2}) defined over~$\mathbb Q$, with
Jacobian~$J$, and suppose that
\begin{itemize}
\item[(i)] $F(X)$ has no root $\theta \in {\mathbb Q}$.\hfil
\item[(ii)] The roots of $F(X)$ cannot be divided into two sets of three
roots, where the sets are either defined over~${\mathbb Q}$ (as wholes)
or defined over a quadratic extension and conjugate over~$\mathbb Q$.
\end{itemize}
If the kernel of~$\mu$ is $2J({\mathbb Q})$ then there does not exist
a $\mathbb Q$-rational divisor class of degree~1 on~$\mathcal{C}$.
In particular, $\mathcal{C}({\mathbb Q}) = \emptyset$.
\end{lem}

We are now in a position to give our example. We take a prime $p \equiv
7~(\hbox{mod }8)$, as well as $a\in {\mathbb Z},
a \not= p,2p$, and let
$\mathcal{C}_{p,a}$ denote
\begin{equation}\label{eq:Cpa}
\mathcal{C}_{p,a} : Y^2 = 2(X^2 + p)(X^2 + 2p)(X^2 + a).
\end{equation}

\begin{ex}\label{ex:KERMUimplyNORDCDO}
The curve $\mathcal{C}_{7,-11} : Y^2 = 2(X^2 + 7)(X^2 + 14)(X^2 - 11)$
has points everywhere locally, but has
no~$\mathbb Q$-rational divisor class of degree~1 and
so~${\mathcal{C}}_{7,-11}({\mathbb Q}) = \emptyset$;
the Jacobian~$J_{7,-11}$ is isogenous over~$\mathbb Q$ to
two elliptic curves defined over~$\mathbb Q$, each of rank~1,
and each with finite Tate-Shafarevich group.
Hence $\mathcal{C}_{7,-11}$ satisfies the conditions
of Proposition~1.3.
\end{ex}
\begin{dem}
\ First note that, for all~${\mathbb Q}_p$  at least one of $\gamma_1,\gamma_2,\gamma_3\in
{\mathbb Q}_p$ exist such that $\gamma_1^2 = 2, \gamma_2^2 = -7$ or $\gamma_3^2=-14$, and
so at least one of $\infty^+, (\gamma_2,0)$ or $(\gamma_3,0) \in
\mathcal{C}_{7,-11}({\mathbb Q}_p)$; a similar statement holds over $\mathbb R$.
Hence~$\mathcal{C}_{7,-11}$ has points everywhere locally. Furthermore, there are maps
$(x,y) \mapsto (x^2,y)$ and $(x,y) \mapsto (1/x^2,y/x^3)$ from $\mathcal{C}_{7,-11}$ to
the elliptic curves
$$Y^2 = 2(X + 7)(X + 14)(X - 11)$$
and $$Y^2 = 2(1 + 7X)(1 + 14X)(1 - 11X),$$ and $J_{7,-11}$ is
isogenous over~$\mathbb Q$ to their product (a special case
of the isogeny in Theorem~14.1.1 of~\cite{prolegom}).

Furthermore, each of these elliptic curves can be shown to have analytic rank~1 (for
example, using Magma~\cite{magma}) and so by Theorem~\ref{thm:koly} they each have rank~1
and finite Tate-Shafarevich group. Since their product is isogenous over~$\mathbb Q$
to~$J_{7,-11}$, it also follows that the rank of $J_{7,-11}({\mathbb Q})$ is~2. Clearly
$J_{7,-11}({\mathbb Q})[2] \cong {\mathbb Z}/2{\mathbb Z} \times {\mathbb Z}/2{\mathbb
Z}$, since points of order~2 in $J_{7,-11}({\mathbb Q})$ correspond to quadratics $q(X)$,
defined over~$\mathbb Q$, with $q(X) | 2(1 + 7X)(1 + 14X)(1 - 11X)$ (see p. 3
of~\cite{prolegom}) and there are three such~$q(X)$. So (since the rank is~2), we have
that $J_{7,-11}({\mathbb Q})/2J_{7,-11}({\mathbb Q}) \cong {\mathbb Z}/2{\mathbb Z} \times
{\mathbb Z}/2{\mathbb Z} \times {\mathbb Z}/2{\mathbb Z} \times {\mathbb Z}/2{\mathbb Z}$.
It only remains to show that ${\mathcal C}_{7,-11}$ has no rational divisor class of
degree~1. After a short search on~$J_{7,-11}({\mathbb Q})$ one finds
$$
\begin{array}{rl}
T_1 &= \{ (\sqrt{-7}, 0), (-\sqrt{-7}, 0)\},\
T_2 = \{ (\sqrt{-14}, 0), (-\sqrt{-14}, 0)\},\\ & \\
D_1 &= \{ (\sqrt{-\frac{23}{2}}, \frac{45}{2} \sqrt{-\frac{23}{2}}),
(-\sqrt{-\frac{23}{2}}, -\frac{45}{2} \sqrt{-\frac{23}{2}}) \},\\ & \\
D_2 &= \{ (\sqrt{-\frac{77}{23}}, \frac{420}{23} \sqrt{-\frac{77}{23}}),
         (-\sqrt{-\frac{77}{23}}, -\frac{420}{23} \sqrt{-\frac{77}{23}})\}.
\end{array}
$$
Applying the map~$\mu$ of~(\ref{eq:mumapgen2})
to~$n_1 T_1 + n_2 T_2 + n_3 D_1 + n_4 D_2$,
for all~16 choices of~$n_i =0,1$, we find that only the
case~$n_1 = n_2 = n_3 = n_4 = 0$ is mapped by~$\mu$ to
the identity and so $T_1,T_2,D_1,D_2$ are independent
in $J_{7,-11}({\mathbb Q})/(\hbox{ker}\,\mu)$. Since
$2J_{7,-11}({\mathbb Q}) \subseteq \hbox{ker}\,\mu$ these must
also be independent in $J_{7,-11}({\mathbb Q})/2J_{7,-11}({\mathbb Q})$;
since also $J_{7,-11}({\mathbb Q})/2J_{7,-11}({\mathbb Q}) \cong
{\mathbb Z}/2{\mathbb Z} \times {\mathbb Z}/2{\mathbb Z} \times
{\mathbb Z}/2{\mathbb Z} \times {\mathbb Z}/2{\mathbb Z}$ it follows that
$T_1,T_2,D_1,D_2$ generate $J_{7,-11}({\mathbb Q})/2J_{7,-11}({\mathbb Q})$
and that $\hbox{ker}\,\mu = 2J_{7,-11}({\mathbb Q})$.
The roots of $2(X^2 + 7)(X^2 + 14)(X^2 - 11)$ clearly satisfy~(i),(ii)
of Lemma~\ref{lem:SMALLKERMUimplyNORDCDO}, from which we
deduce that there does not exist a $\mathbb Q$-rational divisor
class of degree~1 on~$\mathcal{C}_{7,-11}$, as required.
\end{dem}
\par
All curves in the family~(\ref{eq:Cpa}) have points everywhere locally,
by the same argument as given for Example~\ref{ex:KERMUimplyNORDCDO}.
The above example was found by initially searching
in Magma~\cite{magma} the family~(\ref{eq:Cpa}), with
$p\leq 50$ and $-20 \leq a \leq 20$, for cases where
$J_{p,a}$ is isogenous over~$\mathbb Q$ to the product of two
rank~1 elliptic curves, and $\hbox{ker}\,\mu = 2J_{p,a}({\mathbb Q})$.
There were eight such pairs $(p,a)$:
$$
(7,-19),(7,-11),(23,13),(31,-14),
(31,-11),(31,5),(31,13),(47,13).
$$
For all of these pairs, $\mathcal{C}_{p,a}$ has no
$\mathbb Q$-rational divisor class of degree~1, by the same
argument as given for Example~\ref{ex:KERMUimplyNORDCDO},
and so all of these satisfy the conditions of Proposition~1.3.
\par
The author acknowledges the helpful advice
of Nils Bruin and Michael Stoll. He thanks EPSRC for support: grant number EP/F060661/1.


\begin{thebibliography}{9}

\bibitem{bash64} M. I. Bashmakov, On the divisibility of principal homogeneous spaces over Abelian varieties
 (Russian),   {\em Izv. Akad. Nauk SSSR Ser. Mat.}  28  (1964), 661--664.

\bibitem{bash72} \bysame, Cohomology of Abelian varieties over a number field (Russian), {\em Uspehi Mat. Nauk}  27  (1972), 25--66.

\bibitem{bcts} M. Borovoi, J.-L. Colliot-Th\'el\`ene, A. N. Skorobogatov,
The elementary obstruction for homogeneous spaces, {\em Duke Math. J.} 141 (2008),
321--364.

\bibitem{bcdt:mod} C. Breuil, B. Conrad, F. Diamond and R. Taylor.
On the modularity of elliptic curves over ${\mathbb Q}$:
wild 3-adic exercises, {\em J. Amer. Math. Soc.} 14 (2001), 843--939.

\bibitem{bruinstoll:manin} N. Bruin and M. Stoll,
Deciding existence of rational points on curves: an experiment, {\em Experiment.\ Math.}
17 (2008), 181--189.

\bibitem{prolegom} J. W. S. Cassels and E. V. Flynn, {\em Prolegomena to a
    Middlebrow Arithmetic of Curves of Genus 2}, London Math. Soc.
  Lecture Notes Ser., vol. 230, Cambridge University Press, 1996.

\bibitem{es} D. Eriksson, V. Scharaschkin, On the Brauer-Manin obstruction for zero-cycles on curves, {\em Acta Arithmetica} 135 (2008), 99--110.

\bibitem{flynn:manin} E. V. Flynn, The Hasse principle and the Brauer-Manin obstruction
for curves, {\em Manuscripta Math.}  115 (2004), 437--466.

\bibitem{bf} A. Grothendieck, Brief an G. Faltings, reprinted in P. Lochak, L. Schneps (eds.) {\em Geometric Galois actions I},  London Math. Soc. Lecture Note Ser., vol. 242,  Cambridge University Press, Cambridge, 1997,  pp. 49--58; English translation {\em ibid.}, pp. 285--293.

\bibitem{hasza} D. Harari, T. Szamuely, Local-global principles for 1-motives, {\em Duke Math. J.} 143 (2008), 531--557.

\bibitem{jannsen} U. Jannsen, Continuous \'etale cohomology, {\em Math. Ann.} 280 (1988), 207--245.

\bibitem{kala} N. Katz, S. Lang, Finiteness theorems in geometric class field theory, {\em L'Enseign. Math.} 27 (1981), 285--314

\bibitem{koenigsmann} J. Koenigsmann, On the section conjecture in anabelian geometry, {\em J.
reine angew. Math.} 588 (2005), 221--235.

\bibitem{koly:anal} V. A. Kolyvagin, Euler Systems. in P. Cartier et
  al. (eds.) {\em The Grothendieck Festschrift}, volume II,
Progress in Mathematics, vol. 87, Birkh\"auser, Boston, 1990, 435--483.

\bibitem{licht} S. Lichtenbaum, Duality theorems for curves over $p$-adic fields, {\em
Invent. Math.} 7 (1969), 120--136.

\bibitem{magma} The Magma Computational Algebra System,

{\tt http://magma.maths.usyd.edu.au/magma/}

\bibitem{mattuck} A. Mattuck, Abelian varieties over $p$-adic ground fields, {\em Ann. of
Math.} 62 (1955), 92--119.

\bibitem{milne} J. S. Milne, Comparison of the Brauer group with the Tate-\v Safarevi\v c group,  {\em J. Fac. Sci. Univ. Tokyo Sect. IA Math.}  28  (1981), 735--743.

\bibitem{adt} \bysame, {\em Arithmetic Duality Theorems}, second edition, BookSurge, LLC, Charleston, 2006.

\bibitem{poonen} B. Poonen, M. Stoll, The Cassels-Tate pairing on polarized abelian varieties, {\em Ann. of Math.} {\bf 150} (1999), 1109--1149.


\bibitem{sharif} S. Sharif, Curves with prescribed period and index over
local fields, {\it J. Algebra} 314 (2007), 157--167.

\bibitem{skobook} A. N. Skorobogatov, {\em Torsors and rational points}, Cambridge Tracts in Mathematics, vol. 144, Cambridge University Press, Cambridge, 2001.

\bibitem{skoro} \bysame, On the elementary obstruction to the existence of rational points (Russian),  {\em Mat. Zametki}  81  (2007), 112--124;  English translation in  {\em Math. Notes}  81  (2007), 97--107.

\bibitem{ssz} M. Spie\ss, T. Szamuely, On the Albanese map for smooth quasi-projective varieties, {\em Math. Ann.} 325 (2003), 1--17.

\bibitem{stix} J. Stix, On the period-index problem in light of the section conjecture, preprint {\tt arXiv:0802.4125}.

\bibitem{tama} A. Tamagawa, The Grothendieck conjecture for affine curves,
  {\em Compositio Math.}  109  (1997), 135--194.

\bibitem{witt} O. Wittenberg, On Albanese torsors and the elementary obstruction to the existence of 0-cycles of degree 1, {\em Math. Ann.} 340 (2008), 805--838.
\end{thebibliography}
\end{document}